\documentclass[12pt]{article}
\usepackage{amsmath, latexsym, amssymb}

\numberwithin{equation}{section}
\newcommand{\R}{\mathbb R}
\newcommand{\C}{\mathbb C}

\newcommand{\N}{\mathbb N}

\newcommand{\Ta}{T_0}

\newcommand{\T}{T_0'}

\newcommand{\dual}{\langle\tau',\tau\rangle}
\newcommand{\zdual}{\langle b',b\rangle}

\newcommand{\ti}{\mathrm{t}}

\newcommand{\z}{b}
\newcommand{\Z}{B}

\newtheorem{definition}{Definition}[section]

\begin{document}

\begin{center}
{\Large \textbf{Evidence of a Gamma Distribution\\ \vskip .25em  for Prime Powers}} \vskip 4em

{ J. LaChapelle}\\
\vskip 2em
\end{center}

\begin{abstract}
If the prime numbers are pseudo-randomly distributed, then analogy with quantum systems suggests that counting primes might be modeled by a non-homogeneous Poisson process. Consequently, postulating underlying gamma statistics, more-or-less standard heuristic arguments borrowed from quantum mechanics in the context of functional integration allows to derive analytic expressions of several average counting functions associated with prime numbers. The expressions are certain sums of incomplete gamma functions that are closely related to logarithmic-type integral functions --- which in turn are well-known to give the asymptotic dependence of the various counting functions up to error terms. The relatively broad success of  quantum heuristics applied to functional integrals in general along with the excellent numerical accuracy of the analytic expressions for the average counting functions provide strong evidence of a gamma distribution for prime powers.

\end{abstract}

\section{Introduction}
It is no stretch to imagine that the prime numbers possess some underlying random nature. On the other hand, average prime-dependent counting functions are remarkably predictable. So it's natural to wonder if methods and insights afforded by quantum physics might contribute to their understanding. Our aim here is to explore this idea using functional/path integral methods.\footnote{The term``path integral" here refers to the functional integral concept used in quantum physics; not a line integral in mathematics.}

Loosely speaking, the Feynman path integral method in quantum mechanics can be described as the evaluation of the exponentiated action functional over a time-ordered graph. The term `path' refers to the original conceptual context which utilizes the position representation of a quantum state. For position-to-position evolutions, the resulting path integrals can be viewed as  expectations of suitably defined Poisson processes. Indeed, as is well-known, the path integral and the Poisson process are just two different representations of the solution to a differential evolution equation.  Meanwhile, a Poisson process can be viewed as the discrete analog of a gamma process. Hence, the expectation of some Poisson process can be formulated as a kind of `propagator', i.e. integral kernel, based on a gamma functional integral \cite{LA1}. So it is not surprising that the solution of the evolution equation can also be represented as a gamma functional integral.

If one is willing to compare the supposedly (pseudo)random occurrence of prime numbers and their predictable averages to quantum evolution, then analogy suggests that counting prime numbers can be formulated as some integral kernel associated with a gamma functional integral. That is, the very nature of the counting process of a prime event dictates a Poisson process which can be represented by a gamma functional integral. If the functional integral is believed to characterize the state-space of prime numbers as it does for the evolution of quantum systems, then gamma statistics rule. All we need do is determine the relevant, perhaps non-homogeneous,
 scaling factor that parametrizes the process.

But gamma does more. Gamma functional integrals play a key role in quantum systems with constraints \cite{LA1},\cite{LA2}. If we insist that we  count only prime numbers, then that can be viewed as a constraint on an appropriate state space enforced by --- fittingly enough --- a gamma functional integral. It is fortunate that a gamma constraint does not alter the underlying gamma statistical nature of the state space. In statical parlance, the gamma distribution is its own conjugate prior. (If it were otherwise, the heuristics would be inconsistent.) From this perspective, \emph{any}  counting process should exhibit underlying gamma statistics, and we use this as our starting hypothesis.

It may be useful to have a physical picture in mind: Consider  a quantum system of two-state (integer/not-integer) `entities' on the positive reals $\R_+$. The observables of interest are projections onto either of the two possible states.  Observation at a random point via a projector gives integer or not integer. Enumeration of the observations then gives a correspondence between the set of natural numbers $\N_+$ and the lattice $\mathbb{Z}_+$, and this correspondence can be used to characterize an eigenstate by its associated natural number --- thus yielding a model of $\mathbb{Z}_+$ in terms of $\N_+$.  Once observed, an eigenstate persists unless the system is perturbed --- perhaps by changing the ``natural" order-isomorphism between the lattice $\mathbb{Z}_+$ and the set of natural numbers $\N_+$ --- which is ill advised and unlikely but allowed in principle. It turns out that the projector onto integers follows the trivial gamma distribution, and enumeration of the integers is given by the trace of the associated propagator over integer states.

Now consider  a quantum system of two-state (prime/not-prime) `entities' localized on the lattice of positive integers $\mathbb{Z}_+$. Counting `prime events' is postulated to be a constrained dynamical random process.\footnote{It is clear that the distribution of prime \emph{numbers} is not random (although there is evidently some random behavior). What we will find is that the evidence suggests the prime \emph{powers} are random variables. To be precise, the physical model posits that \emph{counting prime powers} is a random process following a constrained gamma distribution. It is remarkable that such an arguably preposterous physical model of the set of prime powers in $\mathbb{Z}_+$ actually appears to work.} As in the case of integers, we use the quantum model on $\R_+$ given a suitable projector. Counting primes corresponds to the expectation of an evolution operator, i.e. integral kernel, generated by the projection onto primes, and the integral kernel can be represented as a gamma functional integral.

The physical model requires a discussion of the functional integral representation of a Poisson process which we present in the next section. Given the functional integral representation, reasonable heuristics lead to analytic expressions for the average prime counting function and other prime-related counting functions that are surprisingly accurate.  There are some twists and turns along the way, however. The gamma distribution encodes information about not only primes but prime powers and prime divisors; and Moebius inversion comes into the picture.
With the analytic representation of the average counting functions in hand, we construct some related explicit counting functions --- some old and some new. In the end, the physical model, accurate average counting functions, and explicit formulae together provide a strong argument in favor of a gamma distribution for prime \emph{powers}.

We utilize the Cartier/DeWitt-Morette scheme \cite{CA/D-W3}-\cite{CA/D-M} for functional integration. A brief summary is given in appendix A that will hopefully suffice as background for $\S$2. We want to emphasize from the beginning that, with the exception of the next section defining gamma and Poisson integrators, our presentation is exploratory and at a formal level of mathematical physics since this seems to be a good way to bridge the physics and mathematics.

\section{Gamma and Poisson integrators}
\begin{definition}
Let $\Ta$ be the space of continuous pointed functions
$\tau:(\mathbb{T}_+,\ti_a)\rightarrow(\C_+,1)$ where $\mathbb{T}_+:=[\ti_a,\ti_b]\subseteq\R_+$ and $\C_+:=\R_+\times i\R$. $T_0$ is an abelian group under point-wise multiplication in the first component and point-wise addition in the second. Let $\beta'$ be a fixed element in the dual group $\T$ of linear characters
$\tau':\Ta\rightarrow\C$ and fix a fiducial $\tau_o\in\Ta$ such that $\langle\beta',\tau_o\rangle=c\in\C_+$.

A lower gamma family of integrators
$\mathcal{D}\gamma_{\alpha,\beta',c}(\tau)$ on $T_0$ is characterized
by
\begin{eqnarray}
&&\Theta_{\alpha,\beta'}(\tau,\tau')= e^{ i\dual-\langle\beta',\tau\rangle}\,\tau^\alpha \notag\\
&&Z_{\alpha,\beta',c}(\tau')
=\frac{\gamma\left(\alpha,c\right)}
{\mathrm{Det}(\beta'-i{\tau'})^\alpha}
\end{eqnarray}
where  $\gamma\left(\alpha,c\right)$ is
the lower incomplete gamma functional given by
\begin{equation}
\gamma\left(\alpha,c\right)
=\Gamma(\alpha)e^{-c}\sum_{n=0}^\infty
\frac{(c)^{\alpha+n}}{\Gamma(\alpha+n+1)}\;,
\end{equation}
and the functional determinant $\mathrm{Det}(\beta'-i{\tau'})^\alpha$ is assumed to be well-defined.

An upper gamma family of integrators
$\mathcal{D}\Gamma_{\alpha,\beta',c}(\tau)$ is defined
similarly where
\begin{equation}
\Gamma\left(\alpha,c\right)
=\Gamma(\alpha)-\gamma\left(\alpha,c\right)
\end{equation}
is the upper incomplete gamma functional.
\end{definition}

Take the gamma integrator and regularize by replacing
$\gamma(\alpha,c)$ with the regularized lower
incomplete gamma function
$P(\alpha,c):=\gamma(\alpha,c)/\Gamma(\alpha)$.
Then restrict to the case $\alpha=n\in\mathbb{N}$, $\beta'=\lambda Id'$ with $\lambda\in\C_+$.
\begin{definition}\label{dirac}
Let $\Ta$ be the space of continuous pointed functions
$\tau:(\mathbb{T}_+,\ti_a)\rightarrow(\C_+,1)$
endowed with a lower gamma family of integrators. Let $\alpha=n\in\mathbb{N}$ and
$\langle \beta',\tau_o\rangle=c$
with $c\in\C_+$. The Poisson integrator family
$\mathcal{D}\pi_{n,\beta',c}(\tau)$ is characterized by
\begin{eqnarray}
&&\Theta_{n,\beta'}(\tau,\tau')=e^{i\dual-\langle \beta',\tau\rangle} \tau^n\notag\\
&&Z_{n, \beta',c}(\tau') = \frac{P\left(n, c\right)}
{\mathrm{Det}\left( \beta'-i{\tau'}\right)^n}\;.
\end{eqnarray}
The Poisson family is defined in terms of the primitive
integrator $\mathcal{D}\tau$ by
\begin{equation}
\mathcal{D}\pi_{n,\beta'}(\tau) :=
e^{-\langle \beta',\tau\rangle}\tau^n\,
\mathcal{D}\tau\;.
\end{equation}
\end{definition}
Note the normalization of the fiducial Poisson integrator
\begin{equation}
\int_{{\Ta}}\;\mathcal{D}\pi_{0,\beta',c}(\tau) =1\;,
\end{equation}
and the rest of the family
\begin{equation}
\int_{{\Ta}}\;\mathcal{D}\pi_{n,\beta',c}(\tau)
=P(n, c)\;.
\end{equation}

For quantum physics applications, $\Re(\tau(\ti))=0$ so that $\tau:(\mathbb{T}_+,\ti_a)\rightarrow(i\R,0)$. In this restricted case $T_0$ is a Banach space over $\R$, and  it is useful to define the `shifted' Poisson integrator by
\begin{equation}
\mathcal{D}\widehat{\pi}_{n,\beta',\tau_o}(\tau) :=
e^{-\langle \beta',(\tau-\tau_0)\rangle}\tau^n\,
\mathcal{D}\tau\;.
\end{equation}
The shifted Poisson integrator is used to define the Poisson expectation of $\beta'$ with respect to the fiducial $\tau_0$;
\begin{equation}
{\langle{\beta'}\rangle}_{\tau_o}:=\int_{{\Ta}}
\;\mathcal{D}\widehat{\pi}_{0,\beta',\tau_o}(\tau)=e^{\langle \beta',\tau_0\rangle}\;.
\end{equation}
In particular, if $\langle \beta',\tau_0\rangle=i\int_{\ti_a}^{\ti_b}\beta(\ti)\;d\ti$, then
\begin{eqnarray}
\langle \beta'\rangle_{\tau_o}=
\sum_{n=0}^\infty \frac{1}{n!}
\int_{\ti_a}^{\ti_b}i\beta'(\ti_1)\cdots\int_{\ti_a}^{\ti_b}i\beta'(\ti_n) \;d\ti_1,\ldots,d\ti_n
\end{eqnarray}
where $\ti_a\leq\ti_1<\cdots<\ti_n\leq\ti_b$. Note that $\frac{\partial}{\partial\ti_b}\langle \beta'\rangle_{\tau_o}=i\beta'(\ti_b)\langle \beta'\rangle_{\tau_o}$, so the Poisson expectation solves a first-order evolution equation.

At the other extreme, if $\Im(\tau(\ti))=0$ then $\tau:(\mathbb{T}_+,\ti_a)\rightarrow(\R_+,1)$. We postulate that a sum over the Poisson integrator family in this case can be used to characterize counting processes. Moreover, a sum over Poisson integrators can be represented as a contour integral of the associated $\alpha$-dependent gamma functional integral; the contour integral accounting for the sum over $n$. Consequently, we contend that a counting process can be represented by \emph{a suitably defined} `$\alpha$-trace' of a gamma functional integral.

As a warm-up exercise for the next section, let's calculate the expected number of events occurring up to some cut-off event $c\in\R_+$  for a suitably defined $\alpha$-trace applied to the simple case of an homogenous process. That is, we take $\beta'=Id'$ in the lower gamma integral and restrict the domain of $\tau$ to some fiducial point, say $\ti_b$. In this case, $T_0$ reduces to $\R_+$ so the functional determinant is well-defined, and we get
\begin{eqnarray}\label{counting}
N(c)
:=\mathrm{tr}_\alpha\int_{T_0}
(-1)^\alpha\;\mathcal{D}\gamma_{\alpha,Id',c}(\tau)&:=&\int_{\mathcal{C}}\frac{\Gamma(1-\alpha)}{2\pi
i}\,(-1)^{\alpha}\gamma(\alpha,c)\,d\alpha\notag\\
&=&\sum_{n=1}^\infty\frac{(-1)^{2n}}{(n-1)!}\,\gamma(n,c)\notag\\
&=&\sum_{n=1}^\infty P(n,c)\notag\\
&=&\Gamma(1,-\log(c))=c
\end{eqnarray}
where the contour encloses the non-negative real axis and $\mathrm{Det}(Id')=1$ by definition. Evidently, the trivial gamma distribution  models the positive integers $\mathbb{Z}_+$, and our defined $\alpha$-trace of the homogenous process counts them. More generally, to count positive integers up to the power of a cut-off we use
 \begin{equation}
N(c^r):=\sum_{n=1}^\infty P(n,c^r)=\Gamma(1,-r\log(c))=c^r\;,\;\;\;\;\;\;r\in\mathbb{N}\;.
\end{equation}

\section{Counting primes}
We will use the functional integrals from the previous section to motivate our conjectures for average prime counting functions. For counting process cut-offs, we restrict to the positive reals so $c\in\R_+\subset\C_+$. To be consistent with the prime number literature, we will denote the cut-off by $\mathsf{x}$ instead of $c$ from here on.
\subsection{Average counting}
According to the physical model analogy, the state space $T_0$ represents countable events (up to some cut-off) that may be non-homogeneously scaled if the system is constrained. In particular, restrict the countable events to prime integers. Our hypothesis is this restriction will result in a bounded non-homogenous rate parameter $|\lambda( \mathsf{x})|<\infty$ which can be imposed by a functional constraint. So the average number of primes up to some cut-off $ \mathsf{x}$ is postulated to go like \cite{LA1}
\begin{eqnarray}
\overline{N_p( \mathsf{x})}
&:=&\mathrm{tr}_\alpha\int_{T_0\times C}(-1)^\alpha\;
\mathcal{D}\gamma_{\alpha,Id', \mathsf{x}}(\tau)
\;\mathcal{D}\gamma_{1,i\lambda(\tau)',\infty}(c)\notag\\
&=&\mathrm{tr}_\alpha\int_{T_0}(-1)^\alpha\;
\mathcal{D}\gamma_{\alpha,Id',\lambda( \mathsf{x})}(\tau)\;.
\end{eqnarray}
Moreover, since prime events start at the second positive integer, the proposed $\alpha$-trace that yields the propagator gets shifted according to
\begin{equation}
 \mathrm{tr}_{\alpha+1}\,\mathrm{F}(\alpha)
:=\frac{1}{2\pi i}\int_{\mathcal{C}_{+1}}\frac{\pi\csc(\pi(\alpha+1))}{\Gamma(\alpha+1)}\,\mathrm{F}(\alpha)\;d\alpha
\end{equation}
where the contour encircles the positive integers.

Consequently, we expect for the average number of primes
\begin{equation}
\overline{N_p( \mathsf{x})}
=\sum_{n=1}^\infty\frac{(-1)^{2n+1}}{n!}\,\gamma(n,\lambda( \mathsf{x}))
=-\sum_{n=1}^\infty\frac{1}{n}\,P(n,\lambda( \mathsf{x}))\;.
\end{equation}
The series converges absolutely
\begin{equation}
\lim_{n\rightarrow\infty} \left|\frac{n!}{(n+1)!}
\frac{\left|\gamma(n+1,\lambda( \mathsf{x}))\right|}
{\left|\gamma(n,\lambda( \mathsf{x}))\right|}\right|
=\lim_{n\rightarrow\infty} \left|\frac{1}{(n+1)}\right|
|\lambda( \mathsf{x})|=0\;,
\end{equation}
and note that
\begin{equation}
\overline{N_p( \mathsf{x}+1)}-\overline{N_p( \mathsf{x})}=
-\sum_{n=1}^\infty\frac{1}{n!}
\int_{\lambda( \mathsf{x})}^{\lambda( \mathsf{x}+1)}e^{-t}\,t^{n}\;dt
\sim\frac{-1}{\lambda( \mathsf{x})}\;.
\end{equation}
Empirical evidence therefore suggests the choice $\lambda( \mathsf{x})=-\log( \mathsf{x})$.

Applying similar heuristics to other related average prime counting functions produces the following list of conjectures:\footnote{It is not clear why the sum of primes and prime entropy get an extra factor of $(-1)^{\alpha+1}$ in their respective functional integral representations. A hand waving argument would say that path reversal in these two cases is accompanied by a minus on the summand $p$ which cancels the minus sign coming from path reversal.}
\begin{itemize}
\item Number of primes
  \begin{equation}
  \overline{N_p( \mathsf{x})}=\sum_{n=1}^\infty\frac{(-1)^{1+0}}{n}\,P(n,-\log( \mathsf{x}))
  \approx\sum_{p=2}^{ \mathsf{x}}\,1
  \end{equation}
\item Sum of Primes
 \begin{equation}
    \overline{\sigma_p({ \mathsf{x}})}
    =\sum_{n=1}^\infty\frac{(-1)^{n+0}}{n}
    \,P(n,-\log({{ \mathsf{x}}}))
    \approx\sum_{p=2}^{{ \mathsf{x}}} p
    \end{equation}
\item Chebyshev function
  \begin{equation}
  \overline{{Ch}({ \mathsf{x}})}=\sum_{n=1}^\infty(-1)^{1+1} P(n+1,-\log({ \mathsf{x}}))\approx \sum_{p=2}^{ \mathsf{x}}\log(p)
  \end{equation}
\item Prime entropy
    \begin{equation}
    \overline{ H({ \mathsf{x}})}:=\sum_{n=1}^\infty(-1)^{n+1}\,P(n+1,-\log({ \mathsf{x}}))\approx \sum_{p=2}^{ \mathsf{x}} p\log(p)
    \end{equation}
\end{itemize}

\subsection{Refined conjectures}
We have that $\int_1^{ \mathsf{x}}|d\gamma(n,-\log(t))/dt|\,dt=\gamma(n,-\log( \mathsf{x}))$. Together with the absolute convergence of $\sum\gamma(n,-\log( \mathsf{x}))/n!$, this implies
\begin{equation}
    \overline{N_{p}( \mathsf{x})}=\mathrm{Ei}(\log( \mathsf{x}))-\log(\log( \mathsf{x}))
    \cong\mathrm{li}( \mathsf{x})\;.
\end{equation}
Consequently the conjecture for the average number of primes is not quite correct: In hind sight, $\overline{N_{p}( \mathsf{x})}$ is more like the Moebius inversion of the average number of primes. Also, notice that the $\alpha$-trace (which is perhaps the most contrived step in arriving at $\overline{N_{p}( \mathsf{x})}$) contributes a different collection of terms if a different choice of contour is made. For example, if the contour is shifted to include the origin, then only the $\log\log( \mathsf{x})$ term survives. Similarly, if the original and the shifted contours are subtracted so that the new contour encloses the first event only, then just the $\mathrm{li}( \mathsf{x})$ piece survives.

Evidently the random variable associated with $\tau$ actually represents primes \emph{and} prime powers. And the $\alpha$-trace seems to embody a choice of inclusion-exclusion based on a particular property of the prime powers relative to $ \mathsf{x}$. It is reasonable, therefore, to view $\overline{N_{p}( \mathsf{x})}$ as an approximation to $J( \mathsf{x})-\omega( \mathsf{x})$ where $J( \mathsf{x})$ is Riemann's counting function
\begin{equation}
J( \mathsf{x}):=\sum_{p^k\leq \mathsf{x}}\,\frac{1}{k}=\sum_{n\leq \mathsf{x}}
\frac{\Lambda(n)}{\log(n)}\;,
\end{equation}
and $\omega( \mathsf{x})$ is the weighted sum of prime power divisors of $ \mathsf{x}$
 \begin{equation}
 \omega( \mathsf{x}):=\sum_{p^k\mid \mathsf{x}}\,\frac{1}{k}=\sum_{n\,\mid\, \mathsf{x}}\frac{\Lambda(n)}{\log(n)}\;.
\end{equation}
Hence,  our initial choice for the $\alpha$-trace subtracts out prime powers that divide $ \mathsf{x}$, and the conjecture should be amended to $\overline{N_{p}( \mathsf{x})}=\overline{J( \mathsf{x})}-\overline{\omega( \mathsf{x})}$ where
\begin{equation}
  \overline{J( \mathsf{x})}:=\mathrm{Ei}\left(\log( \mathsf{x})\right)
  \approx\sum_{n\leq \mathsf{x}}
\frac{\Lambda(n)}{\log(n)}
\end{equation}
and
\begin{equation}
\overline{\omega( \mathsf{x})}:=\log(\log( \mathsf{x}))\approx\sum_{n\,\mid\, \mathsf{x}}\frac{\Lambda(n)}{\log(n)}\;.
\end{equation}

Likewise, it is easy to show that $\overline{Ch({ \mathsf{x}})}= \mathsf{x}-\log( \mathsf{x})$. So, according to our interpretation of $\tau$ and the $\alpha$-trace,
\begin{equation}
\sum_{n=1}^\infty P(n+1,-\log({ \mathsf{x}}))\approx \sum_{n\leq \mathsf{x}}\Lambda(n)-\sum_{n\,\mid\, \mathsf{x}}\Lambda(n)\;.
    \end{equation}
Consequently, the initial guess for the Chebyshev function should be replaced by
\begin{equation}
 \overline{Ch( \mathsf{x})}:=\sum_{n=1}^\infty P(n+1,-\log({ \mathsf{x}}))\approx \sum_{\begin{array}{c}
                 \scriptstyle{n\leq \mathsf{x}} \\
                 \scriptstyle{ n\,\nmid\, \mathsf{x}}
               \end{array}}\,\Lambda(n)\;.
\end{equation}
The amended conjecture for this case is
\begin{equation}\label{second Chebyshev}
\overline{\psi( \mathsf{x})}:= \mathsf{x}
 \;\approx\sum_{n\leq \mathsf{x}}\Lambda(n)
\end{equation}
and
\begin{equation}
 \overline{d_\Lambda( \mathsf{x})}:=\log( \mathsf{x})\;\approx\sum_{n\,\mid\, \mathsf{x}}\Lambda(n)\;.
\end{equation}

Similar considerations applied to the sum of primes and prime entropy lead to the conclusion that their respective conjectures actually represent the sum of weighted prime powers and entropy of prime powers that don't divide $ \mathsf{x}$. Moebius inversion of their conjectured sums yield the average sum of primes $\overline{\sigma_{p}( \mathsf{x})}$ and prime entropy $\overline{ H_{p}({ \mathsf{x}})}$.

To summarize, the refined conjectures for the average counting functions are\footnote{Interpreting the related average sums of divisors is left as an exercise.}
\begin{itemize}
\item Number of primes
  \begin{equation}
  \overline{\pi_1( \mathsf{x})}:=\sum_{m=1}^{\infty}\frac{\mu(m)}{m}\,\overline{J( \mathsf{x}^{1/m})}
\end{equation}
\item Sum of primes
  \begin{equation}
  \overline{\sigma_{p}( \mathsf{x})}:=\sum_{m=1}^{\infty}\frac{\mu(m)}{m}\,\overline{J( \mathsf{x}^{2/m})}
\end{equation}
\item First Chebyshev function
  \begin{equation}
  \overline{\theta({ \mathsf{x}})}=\sum_{m=1}^{\infty}\mu(m)\,\overline{\psi( \mathsf{x}^{1/m})}
  \end{equation}
\item Prime entropy
    \begin{equation}
    \overline{ H_{p}({ \mathsf{x}})}:=\frac{1}{2}\sum_{m=1}^{\infty}\mu(m)\,\overline{\psi( \mathsf{x}^{2/m})}
    \end{equation}
\end{itemize}
 It is worth emphasizing that the right-hand sides can be represented in terms of infinite sums of incomplete gamma functions that, nevertheless, converge fairly rapidly.

\subsection{Exact counting}
To obtain explicit integral representations of exact counting functions, we first Moebius invert the average counting functions listed as the refined conjectures. In other words, take the double Moebius-dual of the original conjectured prime power averages. The resulting sums will include prime powers with no exclusions, and these can be transformed using Perron's formula in the usual way.\footnote{To clean up equations a bit, let's agree to absorb the $(2\pi i)^{-1}$ factor multiplying contour integrals into the measure. Our computations will gloss over issues of analysis since rigorous treatments already exist in the literature. Also, we used \emph{Mathematica} to calculate certain residues, and to visualize the explicit formulae.}

The first two examples are well-known. We repeat them here for completeness and for comparison purposes with the gamma hierarchy discussed in the following subsection.

\subsubsection{number of primes}
For the number of primes, start with $ \overline{\pi_1( \mathsf{x})}$. Its Moebius inversion yields\footnote{Recall that $c\in\C_+$ in general, but in this section we have restricted to the case of $c\equiv \mathsf{x}\in\mathbb{N}_+$. Nevertheless, we will continue to write $\mathrm{Ei}\left(\log( \mathsf{x})\right)$ instead of $\mathrm{li}(\mathsf{x})$ to remind of the more general setting.}
\begin{equation}
 \overline{J( \mathsf{x})}:=\sum_{m=1}^\infty\frac{1}{m}\, \overline{\pi_1( \mathsf{x}^{1/m})}=\mathrm{Ei}\left(\log( \mathsf{x})\right)\;.
\end{equation}
The explicit formula is well-known;
\begin{eqnarray}\label{num of primes}
 J( \mathsf{x})&:=&\lim_{T\rightarrow\infty}\int_{c-iT}^{c+iT}\log(\zeta(s))
 \frac{ \mathsf{x}^{s}}{s}\,ds\;\,\;\;\;\;\;c>1\notag\\
 &=&\lim_{T\rightarrow\infty}\int_{c-iT}^{c+iT}\log(\zeta(s))
 \,d\,\mathrm{Ei}\left(\log( \mathsf{x}^s)\right)\;\,\;\;\;\;\;c>1\notag\\
 &=&-\lim_{T\rightarrow\infty}\int_{c-iT}^{c+iT}\mathrm{Ei}\left(\log( \mathsf{x}^s)\right)
 \,d\log(\zeta(s))\;\,\;\;\;\;\;c>1\notag\\
 &=&\mathrm{Ei}(\log( \mathsf{x}))-\sum_\rho\mathrm{Ei}(\log( \mathsf{x}^\rho))-\log(2)
 -\sum_{k=1}^\infty\mathrm{Ei}(\log( \mathsf{x}^{-2k}))\notag\\
 &=&\sum_{n=2}^{ \mathsf{x}}\frac{\Lambda(n)}{\log(n)}=\sum_{p^k\leq \mathsf{x}}\frac{1}{k}\;,
\end{eqnarray}
and the associated explicit prime counting function is therefore $\pi_1( \mathsf{x})=\sum_{m=1}^\infty\frac{\mu(m)}{m}J( \mathsf{x}^{1/m})$.

The boundary term from the integration by parts leading to the third line vanishes because, for $s=c+it$ with $c>1$,
\begin{equation}
\lim_{t\rightarrow\infty}\left|\log(\zeta(c+it))\right|
\leq\lim_{t\rightarrow\infty}\sum_{p^k}\left|\frac{p^{-k(c+it)}}{k}\right|
\leq\sum_{p^k}\frac{p^{-kc}}{k}=\log(\zeta(c))<\infty\;,
\end{equation}
while
\begin{eqnarray}
\lim_{t\rightarrow\infty}\left|\mathrm{Ei}\left(\log( \mathsf{x}^{(c+it)})\right)\right|
&=&\lim_{t\rightarrow\infty}\left|\frac{ \mathsf{x}^{(c+it)}}{(c+it)\log( \mathsf{x})}
\left(1+O\left(\frac{1}{(c+it)\log( \mathsf{x})}\right)\right) \right|\notag\\
&\leq&\lim_{t\rightarrow\infty}\frac{ \mathsf{x}}{\log( \mathsf{x})}
\left|\frac{1}{c+it}
\left(1+O\left(\frac{1}{c+it}\right)\right) \right|=0\;.\notag\\
\end{eqnarray}

The first and third lines of (\ref{num of primes}) express two complementary viewpoints: Explicit formulae for prime-related summatory functions can be obtained through Perron's formula, or they can be obtained by contour integrals with measure $\log(\zeta(s))$. In the first case, one must find the generating series of the associated summatory function, and in the second case one must find the appropriate integrand. Of course it is possible to transform between the two cases so they are functionally equivalent. But their respective interpretations are quite different, and it is useful to possess both perspectives when constructing new explicit formulae.

Parenthetically, referring to appendix \ref{divisor counting}, the exact number of weighted prime factors can be represented by the infinite sum
\begin{equation}
\omega( \mathsf{x})=\sum_{n=1}^\infty\chi(n, \mathsf{x})\,\frac{\Lambda(n)}{\log(n)}\;,
\end{equation}
which, however, receives essentially no contribution for $n> \mathsf{x}$.

\subsubsection{Chebyshev}
For the sum of log primes, start with the first Chebyshev function. Then, according to (\ref{second Chebyshev}), the average second Chebyshev function is simply $\overline{\psi( \mathsf{x})}= \mathsf{x}$. Furthermore, well-known arguments using Perron's formula (see e.g. \cite{MO}) yield the known explicit formula
\begin{eqnarray}
 \psi( \mathsf{x})&:=& -\lim_{T\rightarrow\infty}\int_{c-iT}^{c+iT}\frac{ \mathsf{x}^{s}}{s}\,d\log(\zeta(s))\;\,\;\;\;\;\;c>1\notag\\
 &=& \mathsf{x}-\sum_\rho\frac{ \mathsf{x}^{\,\rho}}{\rho}-\log(2\pi)
 -\frac{1}{2}\log\left(1-\frac{1}{ \mathsf{x}^{\,2}}\right)\notag\\
 &=&\sum_{n=2}^{ \mathsf{x}}\,\Lambda(n)\;.
 \end{eqnarray}
Meanwhile, using (\ref{div counting}),
\begin{equation}\label{log eq.}
d_\Lambda( \mathsf{x})=\log(\lfloor \mathsf{x}\rfloor)=\sum_{n\,\mid\, \mathsf{x}}\Lambda(n)=\sum_{n=1}^\infty\chi(n, \mathsf{x})\Lambda(n)\;.
\end{equation}
Therefore we can associate an explicit Chebyshev counting function based on the gamma distribution given by
\begin{equation}
 Ch( \mathsf{x})=\psi( \mathsf{x})-\log(\lfloor \mathsf{x}\rfloor)=\sum_{\begin{array}{c}
                 \scriptstyle{n\leq \mathsf{x}} \\
                 \scriptstyle{ n\,\nmid\, \mathsf{x}}
               \end{array}}\,\Lambda(n)\;.
\end{equation}
Moebius inversion gives a series representation of $\Lambda$;
\begin{equation}\label{von Mangoldt}
\Lambda( \mathsf{x})=-\sum_{m=1}^\infty\chi(m, \mathsf{x})\,\mu(m)\log(m)\;.
\end{equation}

\subsubsection{weighted sum of prime powers}
To get the sum of primes we need $\sum_{n\leq \mathsf{x}} n\,\frac{\mu(n)\Lambda(n)}{\log(n)}$. But we don't have a zeta generating function for the associated Dirichlet series. Instead, consider $\sum_{n\leq \mathsf{x}} n\,\frac{\Lambda(n)}{\log(n)}$. Its generating series is
\begin{equation}
\log\left(\zeta(s-1)\right)=\sum_{n=1}^\infty\frac{\Lambda(n)}{\log(n)n^{s-1}}\;.
\end{equation}
So we propose the explicit formula
\begin{eqnarray}
 K( \mathsf{x})&:=&\lim_{T\rightarrow\infty}\int_{c-iT}^{c+iT}\log(\zeta(s-1))
 \,d\,\mathrm{Ei}\left(\log( \mathsf{x}^s)\right)\;\,\;\;\;\;\;c>2\notag\\
 &=&-\lim_{T\rightarrow\infty}\int_{c-iT}^{c+iT}\mathrm{Ei}\left(\log( \mathsf{x}^s)\right)
 \,d\log(\zeta(s-1))\;\,\;\;\;\;\;c>2\notag\\
 &=&\mathrm{Ei}(\log( \mathsf{x}^2))-\sum_\rho\mathrm{Ei}(\log( \mathsf{x}^{1+\rho}))-C-\frac{1}{2}
 -\sum_{k=0}^\infty\mathrm{Ei}(\log( \mathsf{x}^{-2k+1}))\notag\\
 &=&\sum_{n=2}^{ \mathsf{x}} n\,\frac{\Lambda(n)}{\log(n)}=\sum_{p^k\leq \mathsf{x}}\frac{p^k}{k}
\end{eqnarray}
with $C=12\log(A)-1$ where $A$ is the Glaisher-Kinkelin constant. Note however, that the Moebius inversion $\sum_{m=1}^\infty\frac{\mu(m)}{m}K( \mathsf{x}^{1/m})$ is not quite the sum of primes.

\subsubsection{weighted entropy of prime powers}
The prime entropy is given by $-\sum_{n\leq \mathsf{x}} n\,\mu(n)\Lambda(n)$, but again we don't have a zeta generating function. So we  settle for $-\sum_{n\leq \mathsf{x}} n\,\Lambda(n)$ with generating function
\begin{equation}
\frac{\zeta'(s-1)}{\zeta(s-1)}=-\sum_{n=1}^\infty\frac{\Lambda(n)}{n^{s-1}}
\end{equation}
and proposed explicit formula
\begin{eqnarray}
\varepsilon( \mathsf{x})&:=&-\lim_{T\rightarrow\infty}\int_{c-iT}^{c+iT}\frac{ \mathsf{x}^s}{s}
 \,d\log(\zeta(s-1))\;\,\;\;\;\;\;c>2\notag\\
 &=&\frac{1}{2} \mathsf{x}^2-\sum_\rho\frac{ \mathsf{x}^{1+\rho}}{1+\rho}-C
 +\arctan\left(\frac{1}{ \mathsf{x}}\right)\notag\\
 &=&\sum_{n=2}^{ \mathsf{x}} n\,\Lambda(n)=\sum_{p^k\leq \mathsf{x}}\frac{p^k}{k}\log(p^k)\;.
\end{eqnarray}
Here also the Moebius inversion $\sum_{m=1}^\infty\mu(m)\varepsilon( \mathsf{x}^{1/m})$ does not give the entropy of primes.

\subsection{The gamma hierarchy}\label{gamma hierarchy}
As already mentioned, the $\alpha$-trace is probably the most contrived step in our heuristic derivations. Its simple contour and measure are artifacts of fitting the gamma process to counting events of an homogenous process so as to identify the natural numbers with the positive integers. There is no reason to expect counting primes to adhere to the $\alpha$-trace, but inspection of the explicit formulae suggests the kernel based on the simple $\alpha$-trace yields what can be interpreted as a probability. For example, note that the double Moebius-dual of $\overline{N_{p}( \mathsf{x})}$ reduces the infinite sum of lower incomplete gamma functions to the single upper incomplete gamma function $\Gamma(0,-\log( \mathsf{x}))=-\mathrm{Ei}(\log( \mathsf{x}))$. Additionally, hindsight suggests that the relevant measure to use goes like $\log(\zeta(s))$. From this perspective, the explicit formula for counting primes looks something like the Mellin inverse of an expectation;
\begin{equation}\label{1}
 J( \mathsf{x})=\lim_{T\rightarrow\infty}\int_{c-iT}^{c+iT}\Gamma\left(0,-\log( \mathsf{x}^s)\right)
 \,d\log(\zeta(s))\;\,\;\;\;\;\;c>1\;.
\end{equation}
Moreover, sums like $\sum_{n\leq \mathsf{x}}n^r\frac{\Lambda(n)}{\log(n)}$ are governed by the same $\Gamma$ probability density but acquire a shifted zeta measure.

The structure of (\ref{1}) allows to construct an entire hierarchy of explicit formulae that includes all of the previous examples;
\begin{eqnarray}\label{2}
 J( \mathsf{x};r,i)&=&-\lim_{T\rightarrow\infty}\int_{c-iT}^{c+iT}\frac{d^i}{ds^i}\Gamma\left(0,-\log( \mathsf{x}^s)\right)
 \,d\log(\zeta(s-r))\;\,\;\;\;\;\;c>r\notag\\
 &=&-\lim_{T\rightarrow\infty}\int_{c-iT}^{c+iT}\frac{d^i}{ds^i}\overline{J( \mathsf{x}^s)}
 \,d\log(\zeta(s-r))\;\,\;\;\;\;\;c>r
\end{eqnarray}
where $r,i\in\mathbb{N}$. Note that $J( \mathsf{x})=J( \mathsf{x};0,0)$, $\psi( \mathsf{x})= J( \mathsf{x};0,1)$,  $K( \mathsf{x})= J( \mathsf{x};1,0)$, and $\varepsilon( \mathsf{x})=J( \mathsf{x};1,1)$. This hierarchy has an obvious Mellin transform interpretation.

To indicate the pattern for higher derivatives;
\begin{eqnarray}
J( \mathsf{x};0,2)&=&-\lim_{T\rightarrow\infty}\int_{c-iT}^{c+iT}
\left(\frac{ \mathsf{x}^s\log( \mathsf{x})}{s}-\frac{ \mathsf{x}^s}{s^2}\right)
 \,d\log(\zeta(s))\notag\\
&=& \mathsf{x}(\log( \mathsf{x})-1)-\sum_\rho\frac{ \mathsf{x}^\rho(\log( \mathsf{x}^\rho)-1)}{\rho^2}+D
 +\sum_{k=1}^\infty\frac{1+\log( \mathsf{x}^{2k})}{(2k)^2 \mathsf{x}^{2k}}\notag\\
 &=&\sum_{n\leq \mathsf{x}}\,\Lambda(n)\log(n)
\end{eqnarray}
and
\begin{eqnarray}
\overline{ J( \mathsf{x};0,2)}=(-1)^{2+1}\sum_{n=1}^\infty\,P(n+2,-\log( \mathsf{x}))= \mathsf{x}(\log( \mathsf{x})-1)
\approx\sum_{n\leq \mathsf{x}}\,\Lambda(n)\log(n)\notag\\
\end{eqnarray}
where $D=\pi^2/12+(\log(\pi)+\log(2))^2-\log(2\pi)^2-\gamma_0^2-2\gamma_1$ with $\gamma_n$ Stieltjes constants.

The compact form of (\ref{2}) suggests a counterpart on the other side of the Moebius inversion. First note the alternative representations
\begin{equation}
\pi( \mathsf{x})=-\sum_{n\leq \mathsf{x}}\frac{\mu(n)\Lambda(n)}{\log(n)}\;,
\end{equation}
and
\begin{equation}
\theta( \mathsf{x})=-\sum_{n\leq \mathsf{x}}\mu(n)\Lambda(n)\;.
\end{equation}
So define the generating function
\begin{equation}
\log\left(\mathfrak{z}(s)\right):=-\sum_{n=1}^\infty\frac{\mu(n)\Lambda(n)}{\log(n)n^{s}}\;,\;\;\;\;\Re(s)>1
\end{equation}
where the convergence follows by comparison with $\zeta(s)$.
Then conjecture the explicit formulae
\begin{equation}\label{3}
 \pi( \mathsf{x};r,i)=\lim_{T\rightarrow\infty}\int_{c-iT}^{c+iT}\frac{d^i}{ds^i}\Gamma(0,-\log(\overline{\theta( \mathsf{x}^s)})
 \,d\log(\mathfrak{z}(s-r))\;\,\;\;\;\;\;c>r\;.
\end{equation}

Of course this representation is only useful if $\mathfrak{z}(s)$ has zeta-like analytic properties that `filter out' contributions from prime powers --- in which case it would presumably help represent the sum of primes and entropy of primes. Notice that
\begin{eqnarray}
\mathfrak{z}(s)=\prod_p\,e^{-p^{-s}}=\prod_p\sum_{k=0}^\infty\frac{(-1)^k}{\Gamma(k+1)}\frac{1}{p^{ks}}
&=&\prod_p\Gamma(1,p^{-s})\;,\;\;\;\;\Re(s)>1\;.
\end{eqnarray}
So we can't hope to generate $\mathfrak{z}(s)$ as a simple Dirichlet series. Further investigation of $\mathfrak{z}(s)$ is warranted.

\section{Conclusion}
Generally reliable heuristics from quantum physics were used to arrive at various conjectured average counting functions associated with prime numbers. The averages are infinite sums of lower incomplete gamma functions, and their accuracy provides fairly compelling evidence that prime-related summatory functions are governed by a non-homogenous gamma distribution.  We should stress that the heuristics are understood to be motivation and not proof. One could just as well postulate the pertinent infinite sums from the beginning. However, the heuristics encode a probability interpretation that provides a framework in which other summatory functions can be constructed and their associated exact counting functions can be relatively simply represented and understood.

There are at least five directions for further work. First, only a handful of summatory functions were considered, and it would be useful to employ similar methods to other number-theoretic functions. Second, rigorous functional analysis of the proposed explicit formulae for the sum of prime powers and entropy of prime powers is necessary to prove their validity. Third, there are several points that remain puzzling: For example, the gamma process seems to include primes and prime powers, and the $\alpha$-trace seems to encode inclusion-exclusion properties. It would be extremely useful to understand more fully the nature of $\tau^\alpha$ for a non-homogenous gamma process for obvious reasons. Fourth, repeating our analysis for counting of prime events on $\mathbb{Z}_+/n\mathbb{Z}_+$ would likely allow the hypothesis to be compared to known results on prime arithmetic progressions.\footnote{Since the technology to transfer integral kernels from a covering space to its base space already exists in a functional integral context, offhand this would seem like a fairly straightforward exercise.} Finally, the gamma distribution perspective suggests avenues of attack on prime $k$-tuple counting.

\appendix
\section{CDM scheme}\label{appendix 1}
The Cartier/DeWitt-Morette scheme \cite{CA/D-W3}-- \cite{CA/D-M}
defines functional integrals in terms of the data
$(\Z,\Theta,\mathrm{Z},\mathbf{F}(\Z))$.

Here $\Z$ is a separable, possibly infinite dimensional Banach space
with a norm $\|\z\|$ where $\z\in \Z$ is an $L^{2,1}$ map
$\z:[\ti_a,\ti_b]\in\R\rightarrow \C^m$. The dual Banach space $\Z'\ni
\z'$ is a space of linear forms with $\zdual_\Z\in\C$ and
induced norm given by
\begin{equation}
  \|\z'\|=\sup_{\z\neq 0}|\zdual|/\|\z\|\;.\notag
\end{equation}
Assume $\Z'$ is separable. Then $\Z'$ is Polish and consequently
admits complex Borel measures $\mu$.

$\Theta$ and $\mathrm{Z}$ are bounded, $\mu$-integrable functionals
$\Theta:\Z\times \Z'\rightarrow \C$ and $\mathrm{Z}:\Z'\rightarrow\C$. The
functional $\Theta(\z,\cdot)$ can be thought of as the functional
analog of a probability distribution function and $\mathrm{Z}(\cdot)$ the
associated characteristic functional.

The final datum is the space of integrable functionals
$\mathbf{F}(\Z)$ consisting of functionals $\mathrm{F}(\z)$ defined relative
to $\mu$ by
\begin{equation}\label{integrable functions}
\mathrm{F}_{\mu}(\z):=\int_{B'}\Theta(\z,\z')\,d\mu(\z')\;.
\end{equation}
Assuming $\mu\mapsto \mathrm{F}_{\mu}$ is injective, then $\mathbf{F}(\Z)$ is a
Banach space endowed with a norm $\|\mathrm{F}_{\mu}\|$ defined to be the
total variation of $\mu$.

The data are used to characterize an integrator
$\mathcal{D}_{\Theta,Z}\z$ on $\Z$ by
\begin{equation}\label{integrator}
  \int_{B}\Theta(\z,\z')\,\mathcal{D}_{\Theta,Z}\z:=\mathrm{Z}(\z')\,.
\end{equation}
This defines a bounded, linear integral operator $\int_B\mathcal{D}_{\Theta,Z}\z$ on the normed
Banach space of integrable functionals $\mathbf{F}(\Z)$ by
\begin{equation}\label{integral definition}
  \int_{B}\mathrm{F}_{\mu}(\z)\,\mathcal{D}_{\Theta,Z}\z:=
  \int_{B'}\mathrm{Z}(\z')\,d\mu(\z')\
\end{equation}
with
\begin{equation}
\left|\int_B
\mathrm{F}_{\mu}(\z)\,\mathcal{D}_{\Theta,Z}\z\right|\leq\|\mathrm{F}_{\mu}\|\;.
\end{equation}

For non-Banach spaces of pointed maps
$M_a$ where now
$m:[\ti_a,\ti_b]\rightarrow\mathbb{U}\subseteq\mathbb{M}$ with
$m(\ti_a)=m_a$ and $\mathbb{U}\subseteq\mathbb{M}$ an open
neighborhood of an arbitrary Riemannian manifold, CDM uses a parametrization
$P:B\rightarrow M_a$ that allows the integral on
$M_a$ to be defined by
\begin{equation}\label{manifold integral}
\int_{M_a}\mathrm{F}(m)\mathcal{D}m:=\int_{B}\mathrm{F}_\mu(P(b))\mathcal{D}_{\Theta,Z}\z
\end{equation}
in complete analogy with the finite dimensional case.

\section{Counting divisors}\label{divisor counting}
The same strategy for counting gamma events can be used to represent sums over divisors. To count divisors, we propose
\begin{eqnarray}\label{div counting}
\sigma( \mathsf{x})
&:=&\sum_{d=1}^\infty \sum_{n=1}^\infty \lfloor \cos^2\left(\frac{\pi  \mathsf{x}}{d}\right) \rfloor P(n,-\log(d))\notag\\
&=&\sum_{d=1}^\infty \lfloor\cos^2\left(\frac{\pi  \mathsf{x}}{d}\right) \rfloor\,\Gamma(1,-\log(d))\notag\\
&=:&\sum_{d=1}^\infty \chi(d, \mathsf{x})\,\Gamma(1,-\log(d))
\end{eqnarray}
with the definition $\chi(d, \mathsf{x}):=\lfloor\cos^2\left(\frac{\pi  \mathsf{x}}{d}\right) \rfloor$. More generally, we can represent the divisor function by
\begin{equation}
\sigma_{r}( \mathsf{x}):=\sum_{d=1}^\infty \chi(d, \mathsf{x})\,\Gamma(1,-r\log(d))\;,\;\;\;\;\;\;r\in\mathbb{N}\;.
\end{equation}
  We won't prove this, but a bit of hand-waving convinces that the sum over non-divisors $\sum_{d\nmid \mathsf{x}}\chi(d, \mathsf{x})\,\Gamma(1,-r\log(d))$ contributes nothing and the remainder just counts divisors to the appropriate power. Of course the choice of indicator function $\chi(d, \mathsf{x})$ is not unique.

\section{Approximate sum of primes and entropy of primes}\label{sums}
Although the sum of primes and entropy of primes don't have obvious explicit formulae, the probability underpinnings allow to guess fairly good estimates.

It is evident that the gamma process with non-homogenous scaling factor $\log( \mathsf{x})$ characterizes prime powers and not simply primes as initially suspected. So it is reasonable to interpret $\overline{\theta( \mathsf{x})}$ as a kind of `shifted cut-off' that subtracts off the contribution of prime powers \emph{in the mean}. Indeed, it turns out that $\sum_{n=1}^\infty\frac{(-1)^n}{n}P(n,-\log( \overline{\theta({ \mathsf{x}})}))$ and $\sum_{n=1}^\infty(-1)^nP(n+1,-\log( \overline{\theta({ \mathsf{x}})}))$ are very good average values of the sum of primes and the entropy of primes. What is interesting is that simply adding the associated sums over zeta zeros \emph{without the shifted cut-off} yields accurate estimates of the exact counting functions. Specifically, our estimates are
\begin{equation}
\sigma_p( \mathsf{x})\cong\mathrm{Ei}(\log(\overline{\theta( \mathsf{x})}^{\,2}))-\sum_\rho\mathrm{Ei}(\log( \mathsf{x}^{1+\rho}))
+\mathrm{small \; terms}
\end{equation}
and
\begin{equation}
H( \mathsf{x})\cong\frac{1}{2}\overline{\theta( \mathsf{x})}^{\,2}-\sum_\rho\frac{ \mathsf{x}^{1+\rho}}{1+\rho}+\mathrm{small \; terms}\;.
\end{equation}
The reader is encouraged to experiment with the graphs of these estimates to appreciate how the zeta zeros ``pull" the average closer to the exact step functions. The estimates clearly do not converge to the exact functions, but they follow the contours surprisingly well. This lends further support to the probability viewpoint.

\end{document}